\PassOptionsToPackage{unicode}{hyperref}
\PassOptionsToPackage{hyphens}{url}
\documentclass[11pt,a4paper]{article} 
\usepackage{longtable}
\usepackage{array}
\usepackage{caption} 
\usepackage{booktabs} 
\usepackage{amsmath,amssymb}
\usepackage{lmodern}
\usepackage{iftex}
\usepackage{needspace}
\ifPDFTeX
  \usepackage[T1]{fontenc}
  \usepackage[utf8]{inputenc}
  \usepackage{textcomp} 
\else 
  \usepackage{unicode-math}
  \defaultfontfeatures{Scale=MatchLowercase}
  \defaultfontfeatures[\rmfamily]{Ligatures=TeX,Scale=1}
\fi 

\usepackage[style=apa,natbib=true]{biblatex}
\bibliography{bibliography}

\IfFileExists{upquote.sty}{\usepackage{upquote}}{}
\IfFileExists{microtype.sty}{
  \usepackage[]{microtype}
  \UseMicrotypeSet[protrusion]{basicmath} 
}{}
\makeatletter
\@ifundefined{KOMAClassName}{
  \IfFileExists{parskip.sty}{%
    \usepackage{parskip}
  }{
    \setlength{\parindent}{0pt}
    \setlength{\parskip}{6pt plus 2pt minus 1pt}}
}{
  \KOMAoptions{parskip=half}}
\makeatother
\usepackage{xcolor}
\usepackage{longtable,booktabs,array}
\usepackage{calc} 
\usepackage{etoolbox}
\makeatletter
\patchcmd\longtable{\par}{\if@noskipsec\mbox{}\fi\par}{}{}
\makeatother
\IfFileExists{footnotehyper.sty}{\usepackage{footnotehyper}}{\usepackage{footnote}}
\makesavenoteenv{longtable}
\usepackage{graphicx}
\makeatletter
\def\maxwidth{\ifdim\Gin@nat@width>\linewidth\linewidth\else\Gin@nat@width\fi}
\def\maxheight{\ifdim\Gin@nat@height>\textheight\textheight\else\Gin@nat@height\fi}
\makeatother
\setkeys{Gin}{width=\maxwidth,height=\maxheight,keepaspectratio}
\makeatletter
\def\fps@figure{htbp}
\makeatother
\ifLuaTeX
  \usepackage{selnolig}  
\fi
\IfFileExists{bookmark.sty}{\usepackage{bookmark}}{\usepackage{hyperref}}
\IfFileExists{xurl.sty}{\usepackage{xurl}}{} 
\author{}
\date{}

\usepackage[a4paper,margin=25mm]{geometry}

\usepackage{mathptmx} 

\usepackage[T1]{fontenc}

\begin{document}

{\centering

{\fontsize{16}{19}
\selectfont
\textbf{\textsc{Robust Transit Line Planning with Limited Public Data: A Nairobi Matatu Case Study}}}

\vspace{16pt} 

 Mourchid Adegbindin\textsuperscript{*}
 
\textit{Department of Operations Research and Information Engineering\\
Cornell University, Ithaca, NY, USA, 14853, ma757@cornell.edu} 

Brenda Dietrich

\textit{Department of Operations Research and Information Engineering\\
Cornell University, Ithaca, NY, USA, 14853, bld34@cornell.edu}

Jamol Pender

\textit{Department of Operations Research and Information Engineering\\
Cornell University, Ithaca, NY, USA, 14853, jjp274@cornell.edu}

* Corresponding author

\vspace{12pt} 

\emph{Regular paper submitted for presentation at International Conference on Advanced Systems in Public Transport and TransitData 2027}

\emph{June 21-24, 2027, Stockholm, Sweden}

\vspace{12pt} 

}

\textbf{Abstract}:  Informal transit networks pose distinct
planning challenges due to limited demand data, uncertain
travel times, and unscheduled operations. This paper develops
a data-driven robust line-planning framework for Nairobi's
matatu network using public GTFS, WorldPop, and OpenStreetMap
data. A capacitated mixed-integer model selects directional
lines, service tiers, and fleet allocations under budgeted
travel-time uncertainty. Under a common fleet budget,
deterministic optimization serves 39.7\% of synthetic demand,
compared with 39.0\% for robust optimization and 28.6\%
for the constructed reference. Across 500 disruption
scenarios, robust optimization reduces late-headway incidence
from 52.3\% to 6.0\% and aggregate segment overload by 98.1\%.

\textbf{Keywords}: Informal transit; Robust optimization; Transit line planning; Synthetic demand; Nairobi matatus

\_\_\_\_\_\_\_\_\_\_\_\_\_\_\_\_\_\_\_\_\_\_\_\_\_\_\_\_\_\_\_\_\_\_\_\_\_\_\_\_\_\_\_\_\_\_\_\_\_\_\_

\hypertarget{introduction}{%
\section{Introduction}\label{introduction}}
Informal and semi-formal public transport provides essential
urban mobility across the Global South, yet decentralized
operations complicate network planning. Service frequencies
may depend on passenger accumulation, routes lack reliable
published schedules, and the passenger counts and travel-time
records required for conventional planning models are often
unavailable. Nairobi's matatu system illustrates both the
importance of these services and the potential of open transit
data. The Digital Matatus project demonstrated how mobile
mapping and the General Transit Feed Specification (GTFS)
could document Nairobi's semi-formal network
\citep{williams2015}. Comparative evidence further shows
that informal transport networks can achieve high structural
efficiency through self-organization \citep{mittal2024}.
These findings motivate planning methods that retain the
existing network's spatial structure while improving service
allocation and operational reliability.

Transit network design has a substantial optimization
literature. \citet{ceder1986} examined bus network design
through passenger and operator objectives, while
\citet{borndorfer2007} developed a column-generation
approach that integrates line planning and passenger routing.
\citet{bertsimas2021} subsequently demonstrated large-scale,
data-driven transit network redesign using real-world data
from Boston. Operational uncertainty introduces additional
planning challenges. \citet{bertsimas2004} established a
tractable budgeted-uncertainty framework, while
\citet{yan2013} considered stochastic travel times in robust
bus network design. More recently, \citet{guo2024}
investigated robust transit frequency setting under uncertain
passenger demand. Together, these studies provide established
foundations for line planning and uncertainty-aware
optimization.

Their application to informal transit systems raises a
specific data and modeling challenge: how can planners
integrate incomplete route records, estimated spatial demand,
directional passenger itineraries, and uncertain travel times
within operationally constrained network design? This study
addresses the question: \emph{Can data-driven robust
optimization improve transit line planning in informal
networks using limited publicly available data while
accounting for travel-time uncertainty and operational
constraints?}

We develop and evaluate a robust line-planning framework for
Nairobi's matatu network. First, a directed, route-preserving
network is reconstructed from Digital Matatus GTFS records,
and a synthetic origin-destination (OD) matrix is estimated
from WorldPop population data and OpenStreetMap activity
locations. Second, a capacitated mixed-integer linear
program selects among documented directional lines,
determines service tiers, and allocates vehicles subject to
passenger assignment, segment-capacity, and fleet
constraints. Travel-time uncertainty enters through a
budgeted robust counterpart of route cycle-time requirements.
Third, a constructed existing-line reference, a deterministic
design, and a robust design are evaluated under a common
fleet budget and 500 shared travel-time disruption scenarios.

The results quantify the trade-off between passenger coverage
and operational reliability. Deterministic optimization
serves 39.72\% of synthetic demand, compared with 38.99\%
under robust optimization and 28.59\% for the constructed
existing-line reference. Across 500 common travel-time
disruption scenarios, the robust design reduces mean
late-headway incidence from 52.30\% to 6.02\% and aggregate
segment-capacity overload by 98.1\% relative to the
deterministic solution. These improvements require no
additional fleet resources and entail a reduction of only
0.73 percentage points in nominal passenger coverage.

\hypertarget{Data and network representation}{%
\section{Data and network representation}\label{Data_and_network representation}}

\subsection{Network reconstruction}

The network representation derives from the 2019 Digital
Matatus GTFS feed \citep{williams2015}. Route identifiers,
trip records, and ordered stop sequences define the original
directional services. All spatial operations use the projected
coordinate system EPSG:32737. To consolidate nearby stop
records without merging distinct terminals or intersections,
we apply complete-linkage clustering with a maximum
within-cluster diameter of 75~m. This reduces 4,284 GTFS
stops to 3,193 physical stages.

Let $G=(V,A)$ denote the directed physical-stage graph.
Each arc $(u,v)\in A$ represents consecutive stages along
at least one documented GTFS trip. The resulting network
contains 4,283 distinct directed arcs and preserves all
136 routes and 272 directional trip patterns. Each
directional line retains its ordered stage sequence, which
defines the passenger paths and segment loads in the
optimization model. No reverse vehicle-service arcs are
inferred from geographic proximity.

\begin{figure}[t]
    \centering
    \includegraphics[width=0.62\linewidth]{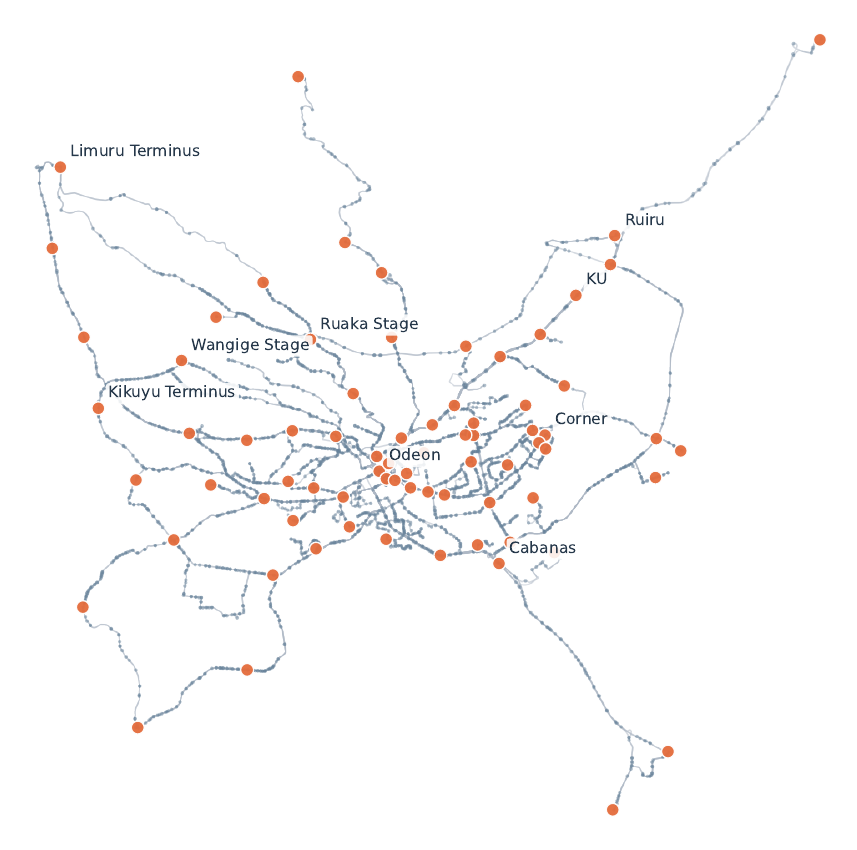}
    \caption{Nairobi matatu network: documented GTFS route
    alignments (grey), consolidated physical stages (blue),
    and strategic hubs (orange). Source: Digital Matatus
    2019 GTFS.}
    \label{fig:network}
\end{figure}

A set of 75 strategic hubs, selected from prominent
terminals and interchange stages, provides the spatial
reference for demand estimation. These hubs constitute a
reduced demand representation; the line-planning model
retains the full physical-stage sequences. Figure~\ref{fig:network}
shows the geographic extent and radial structure of the
documented network.

\subsection{Population and activity data}

The demand model combines two independent spatial inputs.
Residential population derives from the WorldPop Kenya
2026 constrained population-count raster (R2025A v1),
at approximately 100~m resolution \citep{worldpop2025}.
Destination attraction derives from an OpenStreetMap (OSM)
extract of commercial, institutional, educational, and
healthcare facilities \citep{osm2026}.

Each population cell and mapped point of interest (POI)
is assigned to its nearest strategic hub if the Euclidean
distance does not exceed 2.5~km. This procedure assigns
approximately 6.98 million estimated residents and 6,812
of 6,979 processed POIs to the study catchments. The
catchments follow the geographic extent of the matatu
network and do not correspond to Nairobi's administrative
boundary.

POI attraction weights reflect facility categories:
universities receive a weight of 8; colleges, hospitals,
and marketplaces receive 5; schools, clinics, banks, and
offices receive 2; and other mapped activities receive 1.
These weights define relative destination attraction,
rather than measured passenger trip rates.

\subsection{Network access and travel-time representation}

Vehicle-service arcs retain GTFS directionality. Their
nominal travel times derive from straight-line distances
between consecutive stages and an assumed operating speed
of 20~km/h. Separate bidirectional pedestrian links connect
physical stages within 150~m, with an assumed pedestrian
speed of 4.5~km/h and a two-minute transfer allowance.
Passenger itinerary construction permits access to boarding
stages within 750~m of each demand hub.

The 2.5-km catchment radius governs the spatial allocation
of population and POIs, whereas the shorter pedestrian
thresholds govern passenger access and transfers. These
distinctions preserve the separation between geographic
demand estimation and feasible transit service.

The 2026 population estimates and 2019 GTFS routes define
a prospective demand scenario on a documented historical
network. Section~\ref{Synthetic_OD_estimation} describes the resulting
synthetic OD matrix and its sensitivity analysis.

\hypertarget{Synthetic OD estimation}{%
\section{Synthetic OD estimation}\label{Synthetic_OD_estimation}}

We estimate travel demand between the 75 strategic hubs
using the doubly constrained spatial interaction model of
\citet{wilson1967}, with iterative proportional fitting
(IPF) following \citet{morphet1975}. Let $P_i$ denote
the population assigned to origin $i$ and $A_j$ the
weighted POI attraction at destination $j$. The
production and attraction marginals are
\begin{equation}
O_i = Q\frac{P_i}{\sum_k P_k},
\qquad
D_j = Q\frac{A_j}{\sum_k A_k},
\label{eq:marginals}
\end{equation}
where $Q=45{,}000$ trips per hour defines the reference
demand scenario.

Let $c_{ij}$ denote the minimum directed travel impedance
between hubs, inclusive of modeled pedestrian transfers.
The OD matrix satisfies
\begin{equation}
\begin{aligned}
T_{ij} &= a_i b_j O_iD_j
          \exp(-\beta c_{ij}), \qquad T_{ii}=0,\\
\sum_j T_{ij} &= O_i,\qquad
\sum_i T_{ij}=D_j,
\end{aligned}
\label{eq:gravity}
\end{equation}
where $\beta=0.05$~min$^{-1}$ is the reference friction
parameter and $a_i,b_j$ are balancing factors. IPF
enforces both marginals to numerical tolerance.

We test pedestrian-transfer radii
$r\in\{100,150,300,500\}$~m to assess the
sensitivity of the estimated OD distribution.
At 100~m, 148 directed hub pairs lack a modeled
transit-pedestrian path, and 1.84\% of demand
receives a fallback impedance. All hub pairs are
reachable at 150~m and above. Relative to the
150~m case, the alternative OD matrices have
correlations exceeding 0.995, with no more than
1.92\% of total demand reallocated. We retain
150~m as the reference transfer radius, which
achieves full modeled connectivity with the
fewest pedestrian links among the tested
specifications that satisfy this criterion.

\hypertarget{Line and frequency optimization}{%
\section{Line and frequency optimization}\label{Line_and_frequency optimization}}

\subsection{Candidate lines and passenger itineraries}
The candidate set $L$ contains all 272 documented directional
trip patterns over the 3,193 physical stages. For each
positive-demand OD pair $w\in W$, we construct feasible
directed itineraries with at most one transfer, a maximum
hub-to-stage access distance of 750~m, and a pedestrian
transfer radius of 150~m. The procedure produces 7,319
candidate itineraries for 4,218 OD pairs. Each itinerary
retains its ordered line segments, transfer locations,
and generalized access and in-vehicle time. An outside
option accommodates demand without a feasible itinerary
or sufficient service capacity.

\subsection{Mixed-integer formulation}

Let $K=\{\mathrm{low},\mathrm{high}\}$ denote service tiers
with maximum headways $h_k\in\{30,15\}$~min.
For each line $l\in L$, binary variable $y_{lk}$ selects
tier $k$, and integer variable $b_l$ specifies its fleet.
Continuous variables $x_p$ and $u_w$ represent passenger
flow on itinerary $p\in P_w$ and unmet OD demand,
respectively. Let $\tau_p(\Gamma)$ denote itinerary
travel time, including access, transfers, and the
uncertainty adjustment, and let $C_l(\Gamma)$ denote
the corresponding route cycle time.

\Needspace{10\baselineskip}
The model minimizes generalized passenger and operator
cost:
\begin{equation}
\min\
\sum_p \tau_p(\Gamma)x_p
+\frac{1}{2}
\sum_p\sum_{l\in L(p)}\sum_{k\in K}h_k z_{plk}
+180\sum_w u_w
+15\sum_l b_l ,
\label{eq:objective}
\end{equation}
where $L(p)$ is the set of lines used by itinerary $p$
and $z_{plk}=x_p y_{lk}$ for $l\in L(p)$. The product is linearized through
$z_{plk}\le x_p$,
$z_{plk}\le q_w y_{lk}$, and
$z_{plk}\ge x_p-q_w(1-y_{lk})$,
where $w$ is the OD pair of itinerary $p$.
The tier-dependent waiting cost assumes half the
selected maximum headway per boarding.

The principal constraints are
\begin{align}
\sum_{p\in P_w}x_p+u_w &= q_w
&& \forall w\in W, \label{eq:demand}\\
\sum_{k\in K}y_{lk} &\le 1
&& \forall l\in L,\\
b_l &\le B\sum_{k\in K}y_{lk}
&& \forall l\in L,\\
b_l &\ge C_l(\Gamma)
       \sum_{k\in K}\frac{y_{lk}}{h_k}
&& \forall l\in L,\\
\sum_{p:(l,s)\in p}x_p
&\le \frac{60S}{C_l(\Gamma)}b_l
&& \forall(l,s), \label{eq:fleet_capacity}\\
\sum_{p:(l,s)\in p}x_p
&\le 60S\sum_{k\in K}\frac{y_{lk}}{h_k}
&& \forall(l,s), \label{eq:scheduled_capacity}\\
\sum_{l\in L}b_l &\le B.
\label{eq:fleet}
\end{align}

Here $S=33$ seats per vehicle and $B=1{,}000$
vehicles. Passenger flow on each itinerary requires
activation of every line it uses:
$x_p\le q_w\sum_k y_{lk}$ for all $l$ used by $p$.
All flow variables are nonnegative, $y_{lk}$ is
binary, and $b_l$ is nonnegative integer.
The capacity constraint applies separately to every
traversed segment of each directional line.

The selected tier defines the scheduled headway;
additional vehicles provide cycle-time reserve rather
than additional scheduled departures. For an active
line $l$ with selected tier $k$, realized cycle time
$C_l^{(s)}$, and allocated fleet $b_l>0$, the effective
headway in scenario $s$ is
\begin{equation}
h_l^{(s)}=
\max\left\{h_k,\frac{C_l^{(s)}}{b_l}\right\}.
\label{eq:realized_headway}
\end{equation}
This service policy governs the reliability evaluation.

The nominal cycle time equals twice the outbound
running time plus a 10-minute terminal allowance.
This symmetric-return proxy determines fleet
requirements in the absence of observed vehicle duties.
The experiment permits line selection and fleet
reallocation among existing directional patterns;
it does not introduce new corridors.

\subsection{Budgeted travel-time uncertainty}

For each directed stage-arc traversal $e$, travel
time belongs to the uncertainty set
\begin{equation}
t_e=\bar t_e+\hat t_e\xi_e,\quad
0\le\xi_e\le1,\quad
\sum_{e\in E_l}\xi_e\le\Gamma_l,
\label{eq:uncertainty}
\end{equation}
where $\hat t_e=0.25\bar t_e$. The cycle-time proxy
contains two copies of each outbound traversal,
so $|E_l|=n_l$ counts both copies.
We set $\Gamma_l=\gamma n_l$, with
$\gamma=0.25$ for the robust design.

Following \citet{bertsimas2004}, the worst-case
cycle-time increment is
\begin{equation}
\Delta_l(\Gamma_l)=
\sum_{j=1}^{\lfloor\Gamma_l\rfloor}\hat t_{(j)}
+(\Gamma_l-\lfloor\Gamma_l\rfloor)
\hat t_{(\lfloor\Gamma_l\rfloor+1)},
\label{eq:robust}
\end{equation}
where deviations appear in descending order.
Thus $C_l(\Gamma)=\bar C_l+\Delta_l(\Gamma_l)$.
The same budgeted-deviation rule applies separately
to the in-vehicle arcs of each passenger itinerary.
The deterministic model sets $\gamma=0$.

\section{Optimization and reliability results}
\label{sec:results}

\subsection{Line selection and passenger coverage}

We compare three designs under identical synthetic demand
and a fleet limit of 1,000 vehicles: a fixed reference
with all 272 directional lines assigned the 30-minute
service tier, a deterministic line and frequency design,
and a robust design with uncertainty budget $\gamma=0.25$.
The reference permits fleet and passenger-flow allocation
but retains the complete directional line catalog.
All three cases satisfy the same demand-conservation,
service-frequency, and segment-capacity constraints.
Gurobi returns feasible solutions with relative
optimality gaps below 2\%; independent postsolve
constraint checks pass for each case.

\begin{table}[t]
\centering
\small
\caption{Line-planning results under a common
1,000-vehicle fleet budget.}
\label{tab:optimization}
\begin{tabular}{lrrr}
\toprule
Metric & Reference & Deterministic & Robust\\
\midrule
Active directional lines & 272 & 166 & 160\\
Vehicles allocated & 876 & 973 & 972\\
Demand served (pax/h) & 12,866 & 17,875 & 17,546\\
Demand served (\%) & 28.59 & 39.72 & 38.99\\
MIP gap (\%) & 0.00 & 1.98 & 1.36\\
Runtime (s) & 0.7 & 129.7 & 220.3\\
\bottomrule
\end{tabular}
\end{table}

Table~\ref{tab:optimization} shows the value of targeted
service allocation within the documented network.
Deterministic optimization activates 166 directional
lines and increases modeled demand coverage from
28.59\% to 39.72\%, a gain of 11.13 percentage points
over the fixed reference. Robust optimization activates
160 lines and serves 38.99\% of demand, an improvement
of 10.40 percentage points over the reference.
The robust solution serves approximately 328 fewer
passengers per hour than the deterministic solution,
equivalent to 0.73 percentage points of total demand.

Both optimized designs use fewer than 1,000 vehicles.
Their near-identical fleet allocations, 973 and 972
vehicles, indicate that the principal distinction
lies in the distribution of service resources and
the cycle-time requirements imposed on selected lines.

\subsection{Reliability under travel-time disruptions}

We evaluate the three fixed designs across 500 common
disruption scenarios. For each directed physical-stage
arc $e$, a scenario samples
$\epsilon_e\sim\mathcal{U}(0,0.25)$ and sets
$t_e^{(s)}=\bar t_e(1+\epsilon_e)$.
The same arc-level realization applies to every design.
Line selections, vehicle allocations, and passenger
assignments remain fixed throughout the evaluation.

Realized cycle times determine effective headways
according to Eq.~\eqref{eq:realized_headway}.
A headway violation occurs when a line cannot maintain
its selected service tier. Segment overload measures
the sum of positive passenger-flow excess over realized
hourly capacity across individual line segments.

\begin{table}[t]
\centering
\small
\caption{Performance under 500 common travel-time
disruption scenarios. Entries are scenario means.}
\label{tab:reliability}
\begin{tabular}{lrrr}
\toprule
Metric & Reference & Deterministic & Robust\\
\midrule
Late headways (\%) & 30.84 & 52.30 & 6.02\\
Segment overload (pax/h) & 2,398 & 6,459 & 122\\
Generalized time (min) & 49.35 & 44.76 & 43.64\\
\bottomrule
\end{tabular}
\end{table}

The robust design substantially improves operational
reliability (Table~\ref{tab:reliability}). Mean
late-headway incidence falls from 52.30\% under
deterministic optimization to 6.02\%, a relative
reduction of 88.5\%. Aggregate segment overload falls
from 6,459 to 122 passenger-capacity units per hour,
a reduction of 98.1\%. The fixed reference records
30.84\% late headways and 2,398 units of overload.

The deterministic design achieves the highest nominal
coverage but exposes a larger proportion of its active
lines to headway deterioration under travel-time shocks.
The robust design protects scheduled departures through
uncertainty-adjusted cycle-time requirements.
Consequently, its effective service capacity remains
closer to the nominal schedule across the tested
disruptions.

Mean generalized time is 44.76 minutes for the
deterministic design and 43.64 minutes for the robust
design. These averages refer to different assigned
passenger populations and therefore do not establish
a matched-demand travel-time improvement.

\subsection{Operational mechanism}

The optimized designs allocate nearly identical fleet
totals but distribute vehicles differently across the
directional network. Of the 148 lines selected by both
designs, the robust solution assigns more vehicles to
67, fewer to seven, and the same number to 74. Seven
shared lines also shift from the 15-minute to the
30-minute service tier. Overall, the robust design
selects 160 lines and allocates 972 vehicles, compared
with 166 lines and 973 vehicles under deterministic
optimization.

These changes help explain the reliability advantage.
Additional vehicles provide reserve capacity to maintain
scheduled departures when cycle times increase, while
longer scheduled headways reduce the fleet required to
maintain service commitments. Robust optimization thus
adjusts both vehicle allocations and frequency targets
rather than simply increasing fleet size. Across the
tested disruptions, these choices accompany an 88.5\%
reduction in late-headway incidence and a 98.1\%
reduction in aggregate segment overload, at a nominal
coverage cost of 0.73 percentage points.

\subsection{Sensitivity to the uncertainty budget}

We vary the uncertainty-budget fraction over
$\gamma\in\{0,0.10,0.25,0.50,0.75\}$, retaining
the same fleet limit, synthetic demand, and 500
evaluation scenarios. Increasing $\gamma$ generally
reduces nominal demand coverage while improving
operational reliability. At $\gamma=0.10$, coverage
is 39.47\%, with 24.96\% late headways and 1,841
units of aggregate segment overload. At the
reference robust setting $\gamma=0.25$, coverage
is 38.99\%, while late-headway incidence falls
to 6.02\% and overload to 122 units.

At $\gamma=0.50$ and $0.75$, no headway violations
or material segment overload are observed in the
tested scenarios, although nominal coverage falls
to 37.55\% and 37.66\%, respectively. The small
nonmonotonicity in coverage reflects the computed
mixed-integer solutions, which retain nonzero
optimality gaps. These results identify
$\gamma=0.25$ as an intermediate compromise
between nominal coverage and reliability under
the assumed disruption model.

\section{Discussion and conclusion}
\label{sec:conclusion}

This study demonstrates the feasibility of data-driven robust
line planning for informal transit networks with limited public
data. The Nairobi case combines route-preserving network
reconstruction, spatially synthesized passenger demand, and
capacitated optimization under travel-time uncertainty.
Within a common 1,000-vehicle fleet budget, deterministic
optimization increases nominal demand coverage from 28.59\%
to 39.72\%, while robust optimization achieves 38.99\%.
Across 500 common disruption scenarios, the robust design
reduces late-headway incidence by 88.5\% and aggregate
segment-capacity overload by 98.1\% relative to the
deterministic solution.

These results reveal a substantial reliability gain at a
modest cost in nominal coverage. Robust optimization serves
only 0.73 percentage points less demand than deterministic
optimization, yet sharply reduces frequency violations and
capacity overload under travel-time disruptions. Despite
nearly identical aggregate fleet allocations, the two designs
produce markedly different operational outcomes. Explicit
protection against cycle-time uncertainty therefore influences
service reliability through the allocation of vehicles and
service tiers, not simply through total fleet size.

The analysis relies on synthetic OD demand, assumed link
speeds, proximity-based pedestrian connections, and a
symmetric vehicle-cycle approximation. The disruption
scenarios represent controlled travel-time shocks rather
than observed traffic conditions. Future research can
incorporate field-calibrated demand and cycle times,
correlated corridor disruptions, and additional candidate
routes.

The central implication is that nominal fleet capacity
alone does not ensure reliable service. Under a fixed
fleet budget, robust line planning provides a systematic
means to protect operational performance while making
the associated cost in passenger coverage explicit.

\printbibliography
\end{document}